\documentclass[a4paper]{amsart}

\usepackage[matrix,arrow,curve,tips]{xy}
            \SelectTips{cm}{}
\usepackage{amsmath,amssymb,amsthm,amscd}

\newtheorem{dfn}{Definition}[section]
\newtheorem{prop}[dfn]{Proposition}
\newtheorem{theo}[dfn]{Theorem}
\newtheorem{ex}[dfn]{Example}

\newcommand{\CC}{\mathbb{C}}
\newcommand{\cM}{\mathcal{M}}
\newcommand{\cN}{\mathcal{N}}
\newcommand{\ra}{\rightarrow}
\newcommand{\com}{\mathbin{{\scriptstyle \circ }}}
\newcommand{\ten}{\mathbin{\otimes}}
\newcommand{\id}{\mathord{\mathrm{id}}}
\newcommand{\pr}{\mathord{\mathrm{pr}}}
\newcommand{\supp}{\mathord{\mathrm{supp}}}
\newcommand{\uni}{\mathord{\mathrm{uni}}}
\newcommand{\inv}{\mathord{\mathrm{inv}}}
\newcommand{\mlt}{\mathord{\mathrm{mlt}}}
\newcommand{\Cc}{\mathord{\mathcal{C}^{\infty}_{c}}}
\newcommand{\Gc}{\mathord{\Gamma^{\infty}_{c}}}
\newcommand{\EtGPD}{\mathrm{EtGPD}}
\newcommand{\GPD}{\mathrm{GPD}}
\newcommand{\C}{\mathrm{C}}
\newcommand{\D}{\mathrm{D}}
\newcommand{\Rep}{\mathsf{Rep}}

\newcommand{\Mod}{\mathsf{Mod}}
\newcommand{\R}{\mathrm{R_{G}}}
\newcommand{\AdCat}{\mathrm{AdCat}}
\newcommand{\Cat}{\mathrm{Cat}}

\title[]
      {Groupoid representations and modules over the convolution algebras}

\author{J. Kali\v{s}nik}
\address{Institute of Mathematics, Physics and Mechanics,
         University of Ljubljana, Jadranska 19,
         1000 Ljubljana, Slovenia}
\email{jure.kalisnik@fmf.uni-lj.si}

\thanks{This work was supported by
        the Slovenian Ministry of Science}

\begin{document}

\begin{abstract}
The classical Serre-Swan's theorem defines a bijective
correspondence between vector bundles and finitely generated
projective modules over the algebra of continuous functions on
some compact Hausdorff topological space. We extend these results
to obtain a correspondence between the category of representations
of an \'{e}tale Lie groupoid and the category of modules over its
convolution algebra that are of finite type and of constant rank.
Both of these constructions are functorially defined on the Morita
bicategory of \'{e}tale Lie groupoids and the given correspondence
represents a natural equivalence between them.
\end{abstract}

\maketitle

\section{Introduction}

There are many phenomena in different areas of mathematics and
physics that are most naturally described in the language of
smooth manifolds and smooth maps between them. However, some
natural constructions, coming from the theory of foliations or
from Lie group actions, result in slightly more singular spaces
and require a different approach. The Morita category of Lie
groupoids and principal bundles
\cite{Hae58,HilSka,Mo91,MoMr2,Mrc99} provides a natural framework
in which to study many such singular spaces like spaces of leaves
of foliations, spaces of orbits of Lie group actions or for
example orbifolds \cite{Adem,Con94,Mo,MoMr}. A Lie groupoid can be
considered as an atlas for the given singular space. It turns out
that different Lie groupoids represent the same geometric space
precisely when they are Morita equivalent, i.e. when they are
isomorphic in the Morita category of Lie groupoids. For this
reason we are primarily interested in those algebraic invariants
of Lie groupoids that are functorially defined on the Morita
category of Lie groupoids.

Tangent bundles, bundles of higher order tensors, line bundles and
other vector bundles play a central role in the study of smooth
manifolds. The theory of representations of Lie groupoids
naturally extends all these notions to the category of \'{e}tale
Lie groupoids \cite{Adem,MoMr} and shows their close connection to
the theory of Lie groups representations. It encompasses many well
known constructions like equivariant vector bundles
\cite{Ati,Seg}, orbibundles over orbifolds \cite{Adem}, foliated
and transversal vector bundles over spaces of leaves of foliations
\cite{MoMr}, as well as ordinary vector bundles over manifolds or
representations of discrete groups. The construction of the
category of representations of a Lie groupoid is invariant under
the Morita equivalence and thus represents one of the basic
algebraic invariants of Lie groupoids.

The Connes convolution algebra of smooth functions with compact
support \cite{Con94,Mrc07,Ren} on an \'{e}tale Lie groupoid is
another example of such an invariant. Smooth functions with
compact support on a smooth manifold and group algebras of
discrete groups are all special examples of convolution algebras.
Finitely generated projective modules over the algebra of smooth
functions on a compact manifold are closely connected to smooth
vector bundles over that manifold by a theorem due to Serre and
Swan \cite{Serre,Swan}. Representations of a discrete group
$\Gamma$ on complex vector spaces on the other hand correspond
bijectively to modules over the group algebra $\CC[\Gamma]$ of the
group $\Gamma$. One is thus lead to believe that both of these
examples represent a similar type of objects, namely a
correspondence between representations of an \'{e}tale Lie
groupoid and a certain class of modules over its convolution
algebra.

We show in this paper (Theorem \ref{Serre-Swan for groupoids}) how
to extend the classical Serre-Swan's correspondence to the
category of \'{e}tale Lie groupoids. The functor of smooth
sections with compact support defines an equivalence between the
category of representations of an \'{e}tale Lie groupoid $G$ and
the category of modules over the convolution algebra of $G$ which
are of finite type and of constant rank. These modules generalize
the well known finitely generated projective modules over the
algebras of functions and coincide with them if the manifold of
objects of $G$ is compact and connected. For example, such a
correspondence makes it possible to study the category of
orbibundles over a compact orbifold with the tools of
noncommutative geometry, applied to the convolution algebra of the
corresponding orbifold groupoid.

The constructions of the categories of representations and of the
modules of finite type and of constant rank extend to morphisms
from the Morita bicategory of \'{e}tale Lie groupoids to the
$2$-category of additive categories. In this more general context,
the given correspondence can be described (Theorem \ref{Natural
equivalence of morphisms}) as a natural equivalence between these
two morphisms.

\section{Basic definitions and examples}

\subsection{The Morita category of Lie groupoids}

For the convenience of the reader, and to fix the notations, we
begin by summarizing some basic definitions and results concerning
Lie groupoids that will be used throughout this paper. We refer
the reader to one of the books \cite{Mac,MoMr,MoMr2} for a more
detailed exposition and further examples.

A {\em Lie groupoid} over a smooth, second countable, Hausdorff
manifold $M$ is given by a smooth manifold $G$ and a structure of
a category on $G$ with objects $M$, in which all the arrows are
invertible and all the structure maps
\[
\xymatrix{G\times^{s,t}_{M}G \ar[r]^-{\mlt} & G \ar[r]^-{\inv} & G
\ar@<2pt>[r]^{s} \ar@<-2pt>[r]_{t} & M \ar[r]^-{\uni} & G }
\]
are smooth, with the source map $s$ a submersion. We allow the
manifold $G$ to be non-Hausdorff, but we assume that the fibers of
the source map are Hausdorff. We write $G(x,y)=s^{-1}(x)\cap
t^{-1}(y)$ for the set of arrows from $x\in M$ to $y\in M$ and
$G_{x}=G(x,x)$ for the isotropy group of the groupoid $G$ at $x$.
The set $G(x,y)$ is a submanifold of $G$ and the isotropy group
$G_{x}$ is a Lie group. If $g\in G$ is any arrow from $x$ to $y$,
and $g'\in G$ is an arrow from $y$ to $z$, then the product
$g'g=\mlt(g',g)$ is an arrow from $x$ to $z$. The map $\uni$
assigns to each $x\in M$ the identity arrow $1_{x}=\uni(x)$ in
$G$, and we often identify $M$ with its image $\uni(M)$ in $G$.
The map $\inv$ maps each arrow $g\in G$ to its inverse $g^{-1}$.

A Lie groupoid is {\em \'{e}tale} if all of its structure maps are
local diffeomorphisms. A {\em bisection} of an \'{e}tale Lie
groupoid $G$ is an open subset $U$ of $G$ such that both $s|_{U}$
and $t|_{U}$ are injective. To any such bisection $U$ corresponds
a diffeomorphism $\tau_{U}\!:s(U)\ra t(U)$ defined by
$\tau_{U}=t|_{U}\com (s|_{U})^{-1}$. Bisections of the groupoid
$G$ form a basis for the topology on $G$, so in particular they
can be chosen arbitrarily small.

Generalised morphisms \cite{MoMr2,Mrc99} turn out to be the right
notion of a map between Lie groupoids in the representation theory
of Lie groupoids. They are closely connected to groupoid actions
and principal bundles, which we briefly describe.

A smooth {\em left action} of a Lie groupoid $G$ on a smooth
manifold $P$ along a smooth map $\pi\!:P\ra M$ is a smooth map
$\mu\!:G\times^{s,\pi}_{M}P\ra P$, $(g,p)\mapsto g\cdot p$, which
satisfies $\pi(g\cdot p)=t(g)$, $1_{\pi(p)}\cdot p=p$ and
$g'\cdot(g\cdot p)=(g'g)\cdot p$, for all $g',g\in G$ and $p\in P$
with $s(g')=t(g)$ and $s(g)=\pi(p)$. We define right actions of
Lie groupoids on smooth manifolds in a similar way.

Let $G$ be a Lie groupoid over $M$ and let $H$ be a Lie groupoid
over $N$. A {\em principal $H$-bundle over $G$} is a smooth
manifold $P$, equipped with a left action $\mu$ of $G$ along a
smooth submersion $\pi\!:P\ra M$ and a right action $\eta$ of $H$
along a smooth map $\phi\!:P\ra N$, such that (i) $\phi$ is
$G$-invariant, $\pi$ is $H$-invariant and both actions commute:
$\phi(g\cdot p)=\phi(p)$, $\pi(p\cdot h)=\pi(p)$ and
$g\cdot(p\cdot h)=(g\cdot p)\cdot h$ for every $g\in G$, $p\in P$
and $h\in H$ with $s(g)=\pi(p)$ and $\phi(p)=t(h)$, (ii) $\pi:P\ra
M$ is a principal right $H$-bundle:
$(\pr_{1},\eta)\!:P\times^{\phi,t}_{N}H\ra P\times^{\pi,\pi}_{M}P$
is a diffeomorphism.

A map $f\!:P\ra P'$ between principal $H$-bundles $P$ and $P'$
over $G$ is equivariant if it satisfies $\pi'(f(p))=\pi(p)$,
$\phi'(f(p))=\phi(p)$ and $f(g\cdot p\cdot h)=g\cdot f(p)\cdot h$,
for every $g\in G$, $p\in P$ and $h\in H$ with $s(g)=\pi(p)$ and
$\phi(p)=t(h)$. Any such map is automatically a diffeomorphism.
Principal $H$-bundles $P$ and $P'$ over $G$ are isomorphic if
there exists an equivariant diffeomorphism between them. There is
a natural structure of a category (actually a groupoid) on the set
of principal $H$-bundles over $G$ for any two Lie groupoids $G$
and $H$. It has principal $H$-bundles over $G$ as objects and
equivariant diffeomorphisms as morphisms between them.

To any smooth functor $\psi:G\ra H$ there corresponds a principal
$H$-bundle $P(\psi)=M\times_{N}^{\psi,t}H$ over $G$ with the
actions given by the maps $g\cdot(x,h)=(t(g),\psi(g)h)$ for $g\in
G(x,y)$ and $(x,h)\cdot h'=(x,hh')$ for $h,h'\in H$ such that
$s(h)=t(h')$. A principal bundle $P$ is isomorphic to one induced
by a functor if and only if it is trivial, i.e. if there exists a
global smooth section of the bundle $P$.

If $P$ is a principal $H$-bundle over $G$ and if $P'$ is a
principal $K$-bundle over $H$, for another Lie groupoid $K$, one
can define the composition $P\ten_{H}P'$ \cite{MoMr2,Mrc96,Mrc99},
which is a principal $K$-bundle over $G$. It is the quotient of
$P\times^{\phi,\pi'}_{N}P'$ with respect to the diagonal action of
the groupoid $H$. So defined composition is associative only up to
a natural isomorphism.

The {\em Morita category} $\mathsf{GPD}$ of Lie groupoids consists
of Lie groupoids as objects and isomorphism classes of principal
bundles as morphisms between them. The morphisms in $\mathsf{GPD}$
are sometimes referred to as Hilsum-Skandalis maps or generalised
morphisms between Lie groupoids. Two Lie groupoids are Morita
equivalent if they are isomorphic in the category $\mathsf{GPD}$.
The Morita category $\mathsf{EtGPD}$ of \'{e}tale Lie groupoids is
the full subcategory of the category $\mathsf{GPD}$ with \'{e}tale
Lie groupoids as objects. If $G$ and $H$ are \'{e}tale Lie
groupoids and if $P$ is a principal $H$-bundle over $G$, then the
corresponding map $\pi\!:P\ra M$ is automatically a local
diffeomorphism. We will be primarily interested in the Morita
bicategory of Lie groupoids which we describe later on in the
paper.

\subsection{Representations of Lie groupoids}

Lie groupoids admit a twofold interpretation. They can be used to
describe symmetries of fibre bundles in a similar way as Lie
groups are used to study symmetries of topological spaces.
However, the extra transversal part of the structure, which is
encoded in the manifold of objects, makes them a convenient model
for singular geometric spaces such as orbifolds, spaces of leaves
of foliations or spaces of orbits of Lie group actions. The theory
of representations of Lie groupoids provides a unified framework
for the study of vector bundles on such geometric spaces and shows
their intimate connection to representations of Lie groups.

Let $G$ be a Lie groupoid over a smooth manifold $M$ and let $E$
be a smooth complex vector bundle of rank $k$ over $M$. A {\em
representation} of the groupoid $G$ on $E$ is a smooth left action
$\rho\!:G\times_{M}E\to E$, denoted by $\rho(g,v)=g\cdot v$, of
$G$ on $E$ along the bundle projection $p:E\to M$ \cite{MoMr},
such that for any arrow $g\in G(x,y)$ the induced map
$g_{*}\!:E_{x}\to E_{y}$, $v\mapsto g\cdot v$, is a linear
isomorphism.

\begin{ex}\rm
(i) Representations of a {\em unit groupoid} associated to a
smooth manifold $M$ correspond precisely to smooth complex vector
bundles over $M$.

(ii) Let $G$ be a {\em point groupoid} with only one object, i.e.
$G$ is a Lie group $K$. The representation theory of $G$ then
coincides with the representation theory of the Lie group $K$ on
finite dimensional complex vector spaces.

(iii) The {\em pair groupoid} $G=M\times M$ over a smooth manifold
$M$ has both projections as source and target maps and
multiplication defined in a natural way. Every representation of
$G$ on a vector bundle $E$ over $M$ amounts to a natural
identification of all the fibers of $E$ and is thus isomorphic to
a trivial representation.

(iv) Let $G=K\ltimes M$ be the {\em translation groupoid} of a
smooth left action of a Lie group $K$ on a manifold $M$. In this
case the representations of the groupoid $G$ correspond to
$K$-equivariant vector bundles over $M$ \cite{Ati,Seg}.

(v) Any \'{e}tale Lie groupoid $G$ over a manifold $M$ has a
natural representation on the complexified tangent bundle of the
manifold $M$, where the action of any arrow is defined by the
differential of the local diffeomorphism corresponding to some
bisection through that arrow \cite{Adem,MoMr}. The cotangent
bundle and tensor bundles all inherit this natural representation,
so it makes sense to speak of vector fields, differential forms or
riemannian metrics on \'{e}tale Lie groupoids.

(vi) Let $G$ be an {\em orbifold groupoid} (i.e. a proper
\'{e}tale Lie groupoid) over $M$. Such groupoids are used as
models for orbifolds \cite{Adem,Mo,Sat}. Representations of such
groupoids correspond to orbibundles as defined in
\cite{Adem,Ruan,Ruan2}.

(vii) Let $(M,\mathcal{F})$ be a foliated manifold.
Representations of the {\em holonomy groupoid}
$\text{Hol}(M,\mathcal{F})$ are sometimes referred to as
transversal vector bundles, while those of the {\em monodromy
groupoid} $\text{Mon}(M,\mathcal{F})$ are referred to as foliated
vector bundles.
\end{ex}

A {\em morphism} between representations $E$ and $F$ of the
groupoid $G$ over $M$ is a $G$-equivariant morphism $\phi:E\ra F$
of vector bundles. More precisely, $\phi:E\ra F$ is a fiberwise
linear smooth map which commutes with bundle projections and
satisfies $\phi(g\cdot e)=g\cdot\phi(e)$ for all $g\in G$ and all
$e\in E_{s(g)}$. Representations of a Lie groupoid $G$ together
with $G$-equivariant morphisms between them form a category
$\Rep(G)$ of representations of $G$. Categories of (equivariant)
vector bundles over a manifold, categories of orbibundles over an
orbifold or categories of representations of Lie groups are some
examples of categories of representations of Lie groupoids that
arise naturally in various contexts. Direct sums, tensor products,
duals and other operations on vector bundles generalize to
representations of Lie groupoids and turn the category $\Rep(G)$
into an additive category for every groupoid $G$.

Generalised maps between groupoids can be used to pull back
representations in the same sense as vector bundles can be pulled
back along smooth maps. Let $G$ and $H$ be Lie groupoids over $M$
respectively $N$ and let $P$ be a principal $H$-bundle over $G$.
For any representation $E$ of the groupoid $H$ we get the pull
back representation $P\ten_{H}E$ as follows (see \cite{Kal} for
details). The pull back bundle $\phi^{\ast}E=P\times_{N}E$ has a
natural structure of a vector bundle over $P$ with projection onto
the first factor as the projection map. Groupoid $H$ acts
diagonally from the right on the space $\phi^{\ast}E$ along the
fibers of the projection onto $M$ and it is easy to see that the
natural map $P\ten_{H}E=\phi^{\ast}E/H\ra M$ is well defined,
smooth and makes $P\ten_{H}E$ a vector bundle over $M$. Finally,
the action of the groupoid $G$ on the space $P$ induces a
representation of the groupoid $G$ on the bundle $P\ten_{H}E$ by
acting on the first factor.

The construction of pulling back representations along a principal
bundle $P$ extends to a functor from the category of
representations of $H$ to the category of representations of $G$.
Define a representation $\Rep(P)(E)=P\ten_{H}E$ of $G$ for any
representation $E$ of $H$ and a morphism
$\Rep(P)(\phi):P\ten_{H}E\ra P\ten_{H}F$ of representations of $G$
by $\Rep(P)(\phi)(p\ten v)=p\ten\phi(v)$ for any morphism
$\phi:E\ra F$ of representations of the groupoid $H$. We thus
obtain a covariant functor
\[
\Rep(P):\Rep(H)\ra\Rep(G).
\]

One can use an alternative description of the above operation in
the case of trivial bundles, i.e. when the principal bundle comes
from a smooth functor. Suppose that $\psi:G\ra H$ is a smooth
functor between Lie groupoids and let $E$ be a representation of
the groupoid $H$. One defines a representation $\psi^{\ast}E$ of
the groupoid $G$ on the vector bundle $\psi_{0}^{\ast}E$ over $M$
with the action $g\cdot(x,v)=(t(g),\psi(g)v)$ for $g\in G(x,y)$
and $v\in E_{\psi_{0}(x)}$. So defined representation is naturally
isomorphic to the representation $P(\psi)\ten_{H}E$ of $G$ via the
isomorphism $f:P(\psi)\ten_{H}E\ra\psi^{\ast}E$, which sends the
element $(x,h)\ten v$ to the element $(x,hv)$.

\subsection{Convolution algebras and principal bimodules}

A smooth manifold $M$ is closely connected with the commutative
algebra $\Cc(M)$ of smooth functions with compact support on $M$.
By replacing a manifold $M$ with an \'{e}tale Lie groupoid $G$
over $M$ and by defining a proper notion of the convolution
product on the space of functions one obtains the convolution
algebra $\Cc(G)$ of the groupoid $G$. It is in general
noncommutative but it contains the algebra $\Cc(M)$ as a
commutative subalgebra. The above construction naturally extends
to a covariant functor from the Morita category of \'{e}tale Lie
groupoids to the Morita category of algebras, i.e. for each
principal bundle between groupoids one naturally constructs a
bimodule between corresponding algebras. In this subsection we
briefly recall the definition of the convolution algebra
\cite{Con94,CraMoe,Mrc99,Mrc07} assigned to an \'{e}tale, not
necessarily Hausdorff, Lie groupoid and of the principal bimodule
\cite{KalMrc,Mrc99} assigned to a principal bundle.

We first recall the definition of the convolution product on the
vector space $\Cc(G)$ of smooth functions with compact support on
a Hausdorff \'{e}tale Lie groupoid $G$. Define a bilinear
operation on the space $\Cc(G)$ by the formula
\begin{equation} \label{Convolution product}
(ab)(g)=\sum_{g=g'g''}a(g')b(g''),
\end{equation}
for any $a,b\in\Cc(G)$. Equipped with this product the space
$\Cc(G)$ becomes an associative algebra called the {\em
convolution algebra} \cite{Con94} of the \'{e}tale Lie groupoid
$G$.

In the case of a general \'{e}tale Lie groupoid a suitable notion
of a smooth function with compact support on a non-Hausdorff
manifold, as given in \cite{CraMoe}, is needed. Considering that
smooth functions on a Hausdorff manifold $M$ correspond precisely
to the continuous sections of the sheaf of germs of smooth complex
valued functions on $M$ it makes sense to use this alternative
approach to define smooth functions with compact support on an
arbitrary manifold $P$. One first considers the vector space of
all (not-necessarily continuous) sections of the sheaf of germs of
smooth functions on $P$. The trivial extension of any smooth
function with a compact support in a Hausdorff open subset of $P$
naturally represents a section of that sheaf. The vector space
$\Cc(P)$ of smooth functions with compact support on $P$ is then
defined to be the subspace of the space of all sections, generated
by such sections. This definition of the vector space $\Cc(P)$
agrees with the classical one if $P$ is Hausdorff, but in general
there exists no natural multiplicative structure on the space
$\Cc(P)$. The support, i.e. the set where the values of the
section are nontrivial, of any function in $\Cc(P)$ is always a
compact subset of $P$ but not necessarily closed if $P$ is a
non-Hausdorff manifold.

The stalk of the sheaf of germs of smooth functions on $P$ at any
point of $P$ is a commutative algebra with identity, which enables
us to perform pointwise operations such as addition,
multiplication or pullbacks along smooth maps. In particular, for
any non-Hausdorff \'{e}tale Lie groupoid $G$ a formula analogous
to the formula (\ref{Convolution product}) can be used to define
the convolution algebra $\Cc(G)$ of the \'{e}tale Lie groupoid
$G$. We refer the reader to \cite{Mrc07} for details.

\begin{ex}\rm
(i) The convolution product on the space $\Cc(G)$ coincides with
the ordinary pointwise product of functions in $\Cc(M)$ if $G$ is
the unit groupoid associated to a smooth manifold $M$. On the
other hand, if $G$ is the point groupoid of a discrete group
$\Gamma$, it follows $\Cc(G)=\CC[\Gamma]$, where $\CC[\Gamma]$ is
the group algebra of the group $\Gamma$. The subalgebra
$\Cc(M)=\{\lambda\cdot 1_{\Gamma}|\lambda\in\CC\}$ is isomorphic
to complex numbers and is central in $\CC[\Gamma]$.

(ii) The convolution algebra $\Cc(G)$ of the pair groupoid on $n$
points coincides with the algebra of $n\times n$ complex matrices
and contains the subalgebra $\Cc(M)$ of diagonal matrices. This
example shows that $\Cc(M)$ need not lie in the center of the
algebra $\Cc(G)$.

(iii) Let $G$ be the translation groupoid of a smooth action of a
discrete group $\Gamma$ on a manifold $M$. The convolution algebra
$\Cc(G)$ of $G$ is known in the literature \cite{Con94} as the
crossed product algebra $\Gamma\ltimes\Cc(M)$.
\end{ex}

One can naturally extend the convolution algebra construction to a
functor from the Morita category of \'{e}tale Lie groupoids to the
Morita category of algebras. Let $G$ and $H$ be \'{e}tale Lie
groupoids and let $P$ be a principal $H$-bundle over $G$. One can
define convolution actions of the algebras $\Cc(G)$ and $\Cc(H)$
on the vector space $\Cc(P)$ to turn it into a
$\Cc(G)$-$\Cc(H)$-bimodule, which is called the {\em principal
bimodule} associated to the principal $H$-bundle $P$ over $G$. The
functor $\Cc$ from the Morita category of \'{e}tale Lie groupoids
to the Morita category of algebras assigns to every \'{e}tale Lie
groupoid its convolution algebra and to an isomorphism class of a
principal bundle the isomorphism class of the associated principal
bimodule \cite{Mrc99} (see also \cite{KalMrc} for a treatment of
the non-Hausdorff case).

Throughout the rest of the paper we will restrict ourselves to
Hausdorff groupoids for simplicity, although essentially the same
formulas apply in the non-Hausdorff case as well.

\section{Groupoid representations and modules over convolution algebras}

The vector space of sections of a smooth vector bundle $E$ over a
manifold $M$ admits a natural action of the algebra of smooth
functions on $M$. Additional structure of a representation of an
\'{e}tale Lie groupoid $G$ on the bundle $E$ allows a natural
extension of that action to the action of the convolution algebra
of the groupoid $G$. In this section we characterize the modules
over the convolution algebra of the groupoid $G$ that arise in
this fashion from sections of representations of $G$.

Results of this type were first considered by Serre \cite{Serre}
in the the category of algebraic varieties and Swan in the
category of compact Hausdorff topological spaces \cite{Swan}. All
our vector bundles will be assumed to be of globally constant
rank, a condition which is automatically satisfied if the manifold
of objects of the groupoid is connected. Similar results hold
however in the case of vector bundles of globally bounded rank as
well.

\subsection{Module of sections of a representation}

Let $M$ be a smooth, Hausdorff and second countable manifold and
denote by $\Cc(M)$ the algebra of smooth functions with compact
support on $M$. We will denote by $\Rep(M)$ the category of smooth
vector bundles over $M$ since it coincides with the category of
representations of the unit groupoid associated to the manifold
$M$. For any vector bundle $E$ over the manifold $M$, the vector
space $\Gc(E)$ of smooth sections of $E$ with compact support
admits a natural structure of a left $\Cc(M)$-module given by
$(fu)(x)=f(x)u(x)$ for any $f\in\Cc(M)$ and any $u\in\Gc(E)$.
Every morphism $\phi:E\ra F$ of vector bundles over $M$ induces a
homomorphism $\Gc(\phi):\Gc(E)\ra\Gc(F)$ of left $\Cc(M)$-modules
by composing with $\phi$, i.e. $\Gc(\phi)(u)=\phi\circ u$. As a
result we obtain the covariant functor
\[
\Gc=(\Gc)_{M}:\Rep(M)\ra {}_{\phantom{}M}\Mod
\]
from the category of smooth vector bundles over the manifold $M$
to the category of left modules over the commutative algebra
$\Cc(M)$.

Now let $G$ be an \'{e}tale Lie groupoid over $M$. We will
associate an action of the convolution algebra $\Cc(G)$ on the
space of sections $\Gc(E)$ to any representation $E$ of the
groupoid $G$. Define a bilinear map
\[
\Cc(G)\times\Gc(E)\ra\Gc(E)
\]
by the formula
\[
(au)(x)=\sum_{t(g)=x}a(g)(g\cdot u(s(g))),
\]
for $a\in\Cc(G)$ and $u\in\Gc(E)$. Since the function $a\in\Cc(G)$
has a compact support, there are only finitely many $g\in
t^{-1}(x)$ with $a(g)\neq 0$ for each $x\in M$, hence $au$ is a
well defined section of the vector bundle $E$. It remains to be
proven that $au$ belongs to the space $\Gc(E)$ and that the above
map really defines an action of the algebra $\Cc(G)$ on the space
$\Gc(E)$.

Strictly speaking, the above formula only holds for Hausdorff Lie
groupoids. However, by evaluating the sections of the sheaf of
germs of smooth functions, one can use virtually the same formula
for non-Hausdorff groupoids as well.

\begin{prop}
The vector space of sections $\Gc(E)$ has a natural structure of a
left module over the convolution algebra $\Cc(G)$.
\end{prop}

\begin{proof}
First we show that $au$ represents a smooth section of the vector
bundle $E$ for any $a\in\Cc(G)$ and any $u\in\Gc(E)$. We can
decompose any function $a\in\Cc(G)$ as a sum $a=\sum a_{j}$ of
functions, each of which has support contained in some bisection,
so we can assume right from the start that the support of $a$ is
contained in some bisection $U$. Let us denote $W=t(U)$ and
$V=s(U)$. We then have the following commuting diagram
\[
\begin{CD}
E|_{W} @<\mu<<U\times_{V}E @>\text{pr}_{E}>>E|_{V} \\
@VpVV @V\text{pr}_{U}VV @VpVV\\
W @<t|_{U}<<U @>s|_{U}>> V
\end{CD}
\]
of maps of vector bundles. The maps $\mu$ and $\text{pr}_{E}$ are
isomorphisms of vector bundles covering the diffeomorphisms
$t|_{U}:U\ra W$ respectively $s|_{U}:U\ra V$. Let us denote by
$\sigma_{U}=(t|_{U})\com (s|_{U})^{-1}:V\ra W$ the diffeomorphism
corresponding to the bisection $U$, and by
$\tau_{U}=\mu\com(\text{pr}_{E})^{-1}:E|_{V}\ra E|_{W}$ the
corresponding isomorphism of vector bundles. The smooth section
$u|_{V}$ of the bundle $E|_{V}$ gets mapped by the above
isomorphism to the smooth section $u'=\tau_{U}\com
u\com(\sigma_{U})^{-1}$ of the bundle $E|_{W}$. Furthermore, since
$a$ has compact support in $U$, the function $a\com(t|_{U})^{-1}$
has compact support in $W$. One can now express the section $au$
of the bundle $E$ as $(au)(x)=(a\com(t|_{U})^{-1})(x)u'(x)$ to
prove that it is a smooth section of $E$ with compact support in
$W$.

To see that the space of sections $\Gc(E)$ is a module over the
algebra $\Cc(G)$, the equality $a(bu)=(ab)u$ must hold for all
$a,b\in\Cc(G)$ and all $u\in\Gc(E)$. To this effect we compute
\begin{align*}
((ab)u)(x)     & =\sum_{t(g)=x}(ab)(g)(g\cdot u(s(g))) & \\
               & =\sum_{t(g)=x}\Big(\sum_{g=g'g''}a(g')b(g'')\Big)
               (g\cdot u(s(g))) & \\
               & = \sum_{\substack{t(g')=x \\ s(g')=t(g'')}}a(g')b(g'')
               \big((g'g'')\cdot u(s(g''))\big)\;. &
\end{align*}
On the other hand we have
\begin{align*}
(a(bu))(x)      & =\sum_{t(g)=x}a(g)\big(g\cdot\big((bu)(s(g))\big)\big) & \\
               & =\sum_{t(g)=x}a(g)(g\cdot\Big(\sum_{t(g')=s(g)}b(g')
               \big(g'\cdot u(s(g'))\big)\Big) & \\
               & = \sum_{\substack{t(g)=x \\ s(g)=t(g')}}a(g)b(g')
               \big(g\cdot(g'\cdot u(s(g')))\big).\; &
\end{align*}
In the last line we have used the linearity of the map
$g_{*}:E_{x}\ra E_{y}$ for each $g\in G(x,y)$. Since $E$ is a
representation of $G$, the equality $g\cdot(g'\cdot e)=gg'\cdot e$
holds for all pairs of composable arrows $g,g'\in G$ and all $e\in
E_{s(g')}$, thus $(ab)u=a(bu)$.
\end{proof}

By the above procedure we obtain a left module of sections
$\Gc(E)$ over the convolution algebra $\Cc(G)$ for any
representation $E$ of the groupoid $G$. Any morphism $\phi:E\ra F$
of representations of $G$ is in particular a morphism of vector
bundles and therefore produces a homomorphism
$\Gc(\phi):\Gc(E)\ra\Gc(F)$ of $\Cc(M)$-modules. Considering that
the map $\phi$ is fiberwise linear and $G$-equivariant we get the
equalities
\begin{align*}
(\Gc(\phi)(au))(x)  & =\phi\big((au)(x)\big) & \\
               & =\phi\Big(\sum_{t(g)=x}a(g)(g\cdot u(s(g)))\Big) & \\
               & =\sum_{t(g)=x}a(g)(g\cdot\phi(u(s(g)))) & \\
               & =(a\Gc(\phi)(u))(x)\;, &
\end{align*}
for $a\in\Cc(G)$ and $u\in\Gc(E)$. The homomorphism $\Gc(\phi)$ is
therefore a homomorphism of $\Cc(G)$-modules so we have the
covariant functor
\[
\Gc=(\Gc)_{G}:\Rep(G)\ra {}_{\phantom{}G}\mathsf{Mod}
\]
from the category of representations of the groupoid $G$ to the
category of left modules over the convolution algebra $\Cc(G)$ of
$G$.

\subsection{Modules of finite type and of constant rank}

According to the previous subsection we can consider a module over
the convolution algebra of an \'{e}tale Lie groupoid $G$ as a
module of sections of some representation of $G$. However, not
every $\Cc(G)$-module is of this kind and it is not too hard to
find counterexamples. In the following subsection we define and
explain the conditions that characterize the modules of sections
of representations of the groupoid $G$.

Let $M$ be a smooth manifold and denote by $\Cc(M)$ the algebra of
smooth functions with compact support on $M$. There is a standard
bijective correspondence between nontrivial homomorphisms
$\eta:\Cc(M)\ra\CC$ of complex algebras and the points of the
manifold $M$. To any $x\in M$ one associates the evaluation
$\epsilon_{x}:\Cc(M)\ra\CC$ at the point $x$ given by
$\epsilon_{x}(f)=f(x)$ for $f\in\Cc(M)$. Conversely, the kernel of
any nontrivial homomorphism $\eta:\Cc(M)\ra\CC$ is a maximal ideal
of the form $\ker(\eta)=\{f\in\Cc(M)|f(x)=0\}$ for a unique point
$x\in M$, thus $\eta=\epsilon_{x}$. We will use the notation
$I_{x}\Cc(M)=\{f\in\Cc(M)|f(x)=0\}$ for the maximal ideal of
functions that vanish at $x$ and $\Cc(M)(x)=\Cc(M)/I_{x}\Cc(M)$
for the quotient algebra. Evaluation at the point $x$ induces a
canonical isomorphism between the algebra $\Cc(M)(x)$ and the
field of complex numbers.

Now let $G$ be an \'{e}tale Lie groupoid over $M$ and let $\cM$ be
a left $\Cc(G)$-module. It follows that $\cM$ is a $\Cc(M)$-module
as well since $\Cc(M)$ is a subalgebra of the convolution algebra
$\Cc(G)$. The $\Cc(G)$-module $\cM$ is of {\em finite type} if it
is isomorphic, as a $\Cc(M)$-module, to some submodule of the
module $\Cc(M)^{k}$ for some natural number $k$. The
$\Cc(M)$-modules of the form $\Cc(M)^{k}$ correspond precisely to
the modules of sections of trivial vector bundles
$M\times\CC^{k}$, so one can roughly think of modules of finite
type as corresponding to subfamilies of trivial vector bundles.

Now choose an arbitrary point $x\in M$. The $\Cc(M)$-module
$I_{x}\cM=I_{x}\Cc(M)\cdot\cM$ is then a $\Cc(M)$-submodule of
$\cM$ and we denote by $\cM(x)=\cM/I_{x}\cM$ the quotient
$\Cc(M)(x)$-module and consider it as a complex vector space.
Suppose now that the $\Cc(G)$-module $\cM$ is of finite type and
let $\Phi:\cM\ra\Cc(M)^{k}$ be an injective homomorphism of
$\Cc(M)$-modules. For each $x\in M$ we obtain an injective complex
linear map $\Phi(x):\cM(x)\ra\Cc(M)(x)^{k}\cong\CC^{k}$, which
shows that $\cM(x)$ is a finite dimensional complex vector space
for each $x\in M$. We denote by
$\text{rank}_{x}\cM=\dim_{\CC}\cM(x)$ the rank of the module $\cM$
at the point $x\in M$. The $\Cc(G)$-module $\cM$ of finite type is
of {\em constant rank} if the function
$x\mapsto\text{rank}_{x}\cM$ is a constant function. One can
similarly define the notions of modules of locally constant rank
and of modules of globally bounded rank.

Suppose now that $M$ is a smooth, Hausdorff and paracompact
manifold and let $E$ be a vector bundle over $M$. The module
$\Gc(E)$ of sections of the bundle $E$ is a basic example of a
module of finite type and of constant rank. One can see that as
follows. Since $M$ is finite dimensional and paracompact, there
exists a vector bundle $F$ over $M$ such that the bundle $E\oplus
F$ is isomorphic to some trivial vector bundle $M\times\CC^{k}$
over $M$; vector bundles with this property are said to be of
finite type. This property basically follows from the proof of
Lemma $5.9$ in \cite{MilSta}. As a result we obtain the
isomorphism $\Gc(E)\oplus\Gc(F)\cong\Cc(M)^{k}$, i.e. the module
$\Gc(E)$ is of finite type. Furthermore, there is a natural
isomorphism $\Gc(E)(x)\ra E_{x}$ of complex vector spaces for
every $x\in M$, induced by the evaluation at the point $x$, which
shows that the module $\Gc(E)$ is of constant rank.

Our notions of modules of finite type and of constant rank are
closely connected with the classical notions in the Serre-Swan's
theorem. The algebra $\Cc(M)$ is unital precisely when the
manifold $M$ is compact and in this case it makes sense to speak
of free and projective modules over the algebra $\Cc(M)$. Finitely
generated, projective $\Cc(M)$-modules correspond in this case to
the modules of finite type and of constant rank if the manifold
$M$ is connected and to the modules of finite type and of locally
constant rank in general.

\subsection{Equivalence between the categories of representations and of
modules of finite type and of constant rank}

We will denote by $\Mod(G)$ the full subcategory of the category
of left modules over the convolution algebra of the groupoid $G$
consisting of modules of finite type and of constant rank. Since
every module of sections of a representation is such a module, we
have the functor
\[
(\Gc)_{G}:\Rep(G)\ra\Mod(G)
\]
which represents a natural equivalence between the given
categories.

\begin{theo}\label{Serre-Swan for groupoids}
Let $G$ be an \'{e}tale Lie groupoid over a smooth manifold $M$.
The functor $(\Gc)_{G}:\Rep(G)\ra\Mod(G)$ is an equivalence
between the category $\Rep(G)$ of representations of $G$ and the
category $\Mod(G)$ of modules over the convolution algebra
$\Cc(G)$ of the groupoid $G$ which are of finite type and of
constant rank.
\end{theo}

Before we begin with the proof of Theorem \ref{Serre-Swan for
groupoids} we briefly recall the classical version of the
Serre-Swan's theorem in the setting of smooth manifolds and
modules over the algebras of smooth functions with compact
support.

\begin{theo}\label{Serre-Swan classical}
The functor $(\Gc)_{M}:\Rep(M)\ra \Mod(M)$ is an equivalence of
categories for any smooth, Hausdorff and paracompact manifold $M$.
\end{theo}

\begin{proof}
The crucial point in the proof of the theorem is the observation
that every vector bundle over a paracompact manifold is of finite
type, i.e. a subbundle of some trivial bundle. Taking this into
account, basically the same proof as in the Swan's original paper
\cite{Swan} goes through.
\end{proof}

Let $G$ be an \'{e}tale Lie groupoid over $M$. We will prove
Theorem \ref{Serre-Swan for groupoids} by constructing a
quasi-inverse
\[
\R:\Mod(G)\ra\Rep(G)
\]
to the functor $\Gc$ to show that it is an equivalence of
categories. For any $\Cc(G)$-module $\cM$ of finite type and of
constant rank we define a vector bundle $\R(\cM)$ over $M$ as
follows. As a set, the bundle $\R(\cM)$ is defined as a disjoint
union of the spaces $\cM(x)$ for $x\in M$
\[
\R(\cM)=\coprod_{x\in M}\cM(x),
\]
together with a natural projection onto the manifold $M$. To
define a topology and a smooth structure on the space $\R(\cM)$ we
first choose a vector bundle $E$ over $M$ and an isomorphism
$\Phi:\Gc(E)\ra \cM$ of left $\Cc(M)$-modules. Such an isomorphism
exists due to the classical version of Serre-Swan's Theorem
\ref{Serre-Swan classical}. The induced map
$\Phi(x):E_{x}\ra\cM(x)$ is then an isomorphism of complex vector
spaces for each $x$, so we can use the fiberwise linear bijection
\[
\phi=\coprod\Phi(x):E\ra\R(\cM)
\]
to define a structure of a smooth vector bundle over $M$ on the
space $\R(\cM)$. So defined vector bundle structure on the space
$\R(\cM)$ is well defined. Namely, if $E'$ is another vector
bundle over $M$ and if $\Phi':\Gc(E')\ra\cM$ is an isomorphism of
$\Cc(M)$-modules, we obtain the isomorphism
$(\Phi')^{-1}\com\Phi:\Gc(E)\ra\Gc(E')$ of $\Cc(M)$-modules.
Bundles $E$ and $E'$ are therefore isomorphic by Theorem
\ref{Serre-Swan classical} and in turn they define the same vector
bundle structure on the space $\R(\cM)$.

We will next use the extra structure of a $\Cc(G)$-module on the
space $\cM$ to define a representation of the groupoid $G$ on the
vector bundle $\R(\cM)$. Choose any arrow $g\in G(x,y)$ and any
vector $v\in\cM(x)$. We can find an element $m\in\cM$ such that
$v=m(x)$ and a function $a\in\Cc(G)$ with compact support in some
bisection $U$ such that $a(g)=1$. Since $\cM$ is a left
$\Cc(G)$-module, the element $am\in\cM$ is well defined and we
define
\[
g\cdot m(x)=am(y)\in\cM(y).
\]
We will denote by
\[
\mu_{\cM}:G\times_{M}\R(\cM)\ra\R(\cM)
\]
the map defined by the above formula.

\begin{prop}
The map $\mu_{\cM}$ defines a representation of the Lie groupoid
$G$ on the vector bundle $\R(\cM)$ over $M$.
\end{prop}

\begin{proof}
Before we start with the proof of the proposition, we list some
properties of the $\Cc(M)$-modules $\cM(x)$ for $x\in M$ and of
the convolution algebra $\Cc(G)$.

(i) For any $m(x)\in\cM(x)$ and any function $f\in\Cc(M)$ the
action of $f$ on $m(x)$ is just the multiplication by $f(x)$ in
the vector space $\cM(x)$. In particular, $fm(x)=m(x)$ if $f(x)=1$
and $fm(x)=0\in\cM(x)$ if $f(x)=0$.

(ii) Let $a\in\Cc(G)$ be a smooth function with compact support in
the bisection $U$ of the groupoid $G$ and denote $W=t(U)$
respectively $V=s(U)$. For any function $f\in\Cc(M)$ with compact
support in $V$ we have the formula $af=(f\com(\sigma_{U})^{-1})a$,
where $\sigma_{U}:V\ra W$ is the diffeomorphism corresponding to
the bisection $U$ and $f\com(\sigma_{U})^{-1}\in\Cc(M)$ is a
smooth function with compact support in $W$. Moreover, $af=a$ for
every function $f$ which is identically equal to $1$ on the set
$s(\supp(a))\subset V$.

Suppose now that $g\in G(x,y)$ is an arrow from $x$ to $y$ and let
$v\in\cM(x)$ be an arbitrary vector. We will first show that the
element $\mu_{\cM}(g,v)$ is well defined and independent of
various choices of representatives for the elements $g$ and $v$.
To this effect choose an arbitrary function $a\in\Cc(G)$ with
support in a bisection $U$ such that $a(g)=1$ and let $m,m'\in\cM$
both satisfy $m(x)=m'(x)=v$. We then have $(m-m')(x)=0$, so we can
find a function $f\in I_{x}\Cc(M)$ and an element $m''\in\cM$ with
$m-m'=fm''$. Now choose a smooth function $f'\in\Cc(M)$ with
compact support in $V=s(U)$, such that $af'=a$ as in (ii). The
support of the function $f'f$ is then contained in $V$, so by (i)
it follows
\[
a(m-m')(y)=afm''(y)=af'fm''(y)=((f'f)\com(\sigma_{U})^{-1})am''(y)=0
\]
since $((f'f)\com(\sigma_{U})^{-1})(y)=(f'f)(x)=0$. We obtain the
equality $am(y)=am'(y)$ which proves the independence of
$\mu_{\cM}(g,v)$ of the choice of the representative for the
element $v$. Next we prove a similar statement for the choice of
the representative for the arrow $g$. Suppose $a'\in\Cc(G)$ is
another function with $a'(g)=1$ and with support in a bisection
$U'$. One can then find functions $f,f'\in\Cc(M)$ with
$f(y)=f'(y)=1$, such that $fa=f'a'$ is a function with compact
support in the bisection $U\cap U'$. The equalities
\[
am(y)=fam(y)=f'a'm(y)=a'm(y)
\]
then show that the element $\mu_{\cM}(g,v)\in\R(\cM)$ is well
defined.

We next show that $\mu_{\cM}$ defines an action of the groupoid
$G$ on the vector bundle $\R(\cM)$ over $M$. This will be true if
we show that the equalities
\[
(g'g)\cdot v=g'\cdot(g\cdot
v)\hspace{10pt}\text{and}\hspace{10pt}1_{x}\cdot v=v
\]
hold for all arrows $g\in G(x,y)$, $g'\in G(y,z)$ and all
$v\in\cM(x)$. If $a\in\Cc(G)$ with $a(g)=1$ has support in a
bisection $U$ and if $a'\in\Cc(G)$ with $a'(g')=1$ has support in
a bisection $U'$, then $a'a$ has support in the bisection
$U'\times_{M}^{s,t}U$ and $(a'a)(g'g)=1$. We then have
\[
(g'g)\cdot v=a'am(z)=g'\cdot am(y)=g'\cdot(g\cdot v).
\]
where $m\in\cM$ is an arbitrary element such that $m(x)=v$. To
prove the second claim, it is enough to observe that any function
$f\in\Cc(M)\subset\Cc(G)$ with $f(x)=1$ represents the identity
arrow $1_{x}$ at the point $x\in M$. The result then follows
directly from the definition of the action of the algebra $\Cc(M)$
on $\cM(x)$.

It remains to be proven that
$\mu_{\cM}:G\times_{M}\R(\cM)\ra\R(\cM)$ is a smooth map. Choose
any element $(g,v)\in G\times_{M}\R(\cM)$ and a function
$a\in\Cc(G)$ with compact support in a bisection $U'$ such that
$a|_{U}\equiv 1$ for some small neighbourhood $U\subset U'$ of the
arrow $g$. There exist elements $m^{1},\ldots,m^{k}\in\cM$ with
the property that the vectors $\{m^{1}(x),\ldots,m^{k}(x)\}$ form
a basis of the vector space $\cM(x)$ for all $x\in s(U)$. As a
result we obtain smooth functions
$\lambda_{1},\ldots,\lambda_{k}:\R(\cM)|_{s(U)}\ra\CC$,
implicitely defined by the formula
$w=\sum_{i=1}^{k}\lambda_{i}(w)m^{i}(x)$ for any $w\in\cM(x)$
where $x\in s(U)$. Locally, on a neighbourhood
$U\times_{M}\R(\cM)$ of the point $(g,v)$, we have
\[
\mu_{\cM}(h,w)=\sum_{i=1}^{k}\lambda_{i}(w)am^{i}(t(h)),
\]
where $am^{i}$ are smooth sections of the vector bundle $\R(\cM)$.
This concludes the proof of the proposition.
\end{proof}

Now choose left $\Cc(G)$-modules of finite type and of constant
rank $\cM$ and $\cN$ and let $\Phi:\cM\ra\cN$ be a homomorphism
between them. For each $x\in M$ we have the induced linear maps
$\Phi(x):\cM(x)\ra\cN(x)$ that define a fiberwise linear map
\[
\R(\Phi)=\coprod_{x\in M}\Phi(x):\R(\cM)\ra\R(\cN).
\]
The map $\R(\Phi)$ is a smooth morphism of vector bundles, since
it transforms smooth sections of the bundle $\R(\cM)$ to smooth
sections of the bundle $\R(\cN)$. We claim that $\R(\Phi)$ defines
a $G$-equivariant morphism of representations of the groupoid $G$
on $\R(\cM)$ respectively $\R(\cN)$ as described above. To see
this, choose an arrow $g\in G(x,y)$ and a vector
$v\in\R(\cM)_{x}=\cM(x)$. Let $a\in\Cc(G)$ be a function with
support in a bisection $U$ such that $a(g)=1$ and suppose
$m\in\cM$ satisfies $m(x)=v$. The statement then follows from the
equalities
\[
\R(\Phi)(g\cdot
v)=\Phi(y)(am(y))=\Phi(am)(y)=a\Phi(m)(y)=g\cdot\R(\Phi)(v).
\]
The functoriality of the assignment $\Phi\mapsto\Phi(x)$ for each
$x\in M$ extends to the functoriality of the map $\R$, so we have
the functor
\[
\R:\Mod(G)\ra\Rep(G).
\]

\begin{proof}[Proof of Theorem \ref{Serre-Swan for groupoids}]
We will prove the theorem by showing that the functors
$\R:\Mod(G)\ra\Rep(G)$ and $\Gc:\Rep(G)\ra\Mod(G)$ represent
mutual inverses.

We can naturally identify the $\Cc(G)$-module $\cM$ with the space
of sections $\Gc(\R(\cM))$ of the representation $\R(\cM)$ of $G$
by assigning a section $x\mapsto m(x)$ of the bundle $\R(\cM)$ to
the element $m\in\cM$. Denote by
$\epsilon_{\cM}:\cM\ra\Gc(\R(\cM))$ the corresponding isomorphism
of $\Cc(G)$-modules and let
\[
\epsilon:\text{Id}_{\Mod(G)}\Rightarrow\Gc\circ\R
\]
be the corresponding natural equivalence of functors.

Let $E$ be a representation of the groupoid $G$. There is a
natural isomorphism $\Gc(E)(x)\ra E_{x}$ of complex vector spaces
for every $x\in M$, induced by the evaluation at the point $x$,
which induces an isomorphism $\eta_{E}:\R(\Gc(E))\ra E$ of
representations of the groupoid $G$. The natural equivalence of
functors
\[
\eta:\R\circ\Gc\Rightarrow\text{Id}_{\Rep(G)}
\]
together with the equivalence $\epsilon$ shows that $\Gc$ is an
equivalence of categories.
\end{proof}

\section{Equivalence of morphisms between bicategories}

According to the previous sections we can associate to any
\'{e}tale Lie groupoid the additive categories of representations
and of modules of finite type and of constant rank. We would like
to extend those constructions to define functors $\Rep$ and $\Mod$
from the category of \'{e}tale Lie groupoids and principal bundles
to the category of additive categories and functors such that
their values at the groupoid $G$ are $\Rep(G)$ respectively
$\Mod(G)$. Furthermore, the family of functors $(\Gc)_{G}$ from
Theorem \ref{Serre-Swan for groupoids} should represent a natural
equivalence between these two functors.

If one wants to work in the framework of categories, it is
necessary to work with isomorphism classes of representations and
of modules and forget the extra categorical structure on these
objects. A more convenient way of describing our statement uses
the language of bicategories and morphisms between them which we
now briefly summarize.

\subsection{Bicategories and Morphisms}

We first recall the definition and some examples of bicategories
and morphisms between them as given in \cite{Ben,Lein} (see also
\cite{Lein2} for a very concise treatment). A {\em bicategory}
$\C$ consists of
\begin{enumerate}
\item A collection $\C_{0}$ of objects of $\C$.

\item Categories $\C(A,B)$ for each pair $A,B$ of objects of $\C$.
The objects and arrows of the categories $\C(A,B)$ are called
$1$-cells respectively $2$-cells of the bicategory $\C$.

\item For all $A,B,C\in\C_{0}$ there exist bifunctors
$\C(A,B)\times\C(B,C)\ra\C(A,C)$ which define composition of
$1$-cells of $\C$ and horizontal composition of $2$-cells of $\C$.
Furthermore, for each $A\in\C_{0}$ there is the identity $1$-cell
$1_{A}\in\C(A,A)$.
\end{enumerate}
Composition of $1$-cells in a bicategory $\C$ resembles the usual
composition of arrows in ordinary categories, although it is in
general associative only up to natural associativity coherence
isomorphisms, which are given as a part of the structure of $\C$.
Similarly, the identity $1$-cells $1_{A}$ for $A\in\C_{0}$ act as
units only up to predefined unit coherence isomorphisms. All these
coherence isomorphisms have to satisfy some natural pentagon and
triangle axioms (see \cite{Lein2} for details) as to insure that
all the diagrams constructed out of coherence isomorphisms
commute. Apart from the horizontal composition of $2$-cells in $C$
there exists a vertical composition of $2$-cells, defined by the
composition in the categories $\C(A,B)$. Strict bicategories or
$2$-categories are bicategories in which the composition of
$1$-cells is strictly associative with strict identities.

\begin{ex}\rm
(i) The bicategory $\Cat$ consists of categories as objects,
functors between them as $1$-cells and natural transformations
between functors as $2$-cells. For any category $A\in\Cat_{0}$ the
$1$-cell $1_{A}$ is represented by the identity functor on $A$.
All the coherence isomorphisms are identities so $\Cat$ is in fact
a $2$-category. We will denote by $\AdCat$ the full sub
$2$-category consisting of additive categories.

(ii) Next we describe the Morita bicategory $\GPD$ of Lie
groupoids. For any two Lie groupoids $G$ and $H$ there exists a
category $\GPD(G,H)$ with principal $H$-bundles over $G$ as
objects and equivariant diffeomorphisms as morphisms between them.
Composition of $1$-cells $P$ and $Q$ is defined by the tensor
product construction $P\ten_{H}Q$ together with the natural
associtativity coherence isomorphisms
\[
a_{P,Q,R}:(P\ten_{H}Q)\ten_{K}R\ra P\ten_{H}(Q\ten_{K}R)
\]
for $P\in\GPD(G,H)$, $Q\in\GPD(H,K)$ and $R\in\GPD(K,L)$. For any
$G\in\GPD_{0}$ the identity $1$-cell $1_{G}$ is simply the
groupoid $G$ viewed as a principal $G$-bundle over $G$. The unit
coherence isomorphisms
\[
l_{P}:G\ten_{G}P\ra
P\hspace{10pt}\text{and}\hspace{10pt}r_{P}:P\ten_{H}H\ra P
\]
are induced by the action maps of the groupoids on the principal
bundles. If $f:P\ra P'$ and $g:Q\ra Q'$ are equivariant
diffeomorphisms, their horizontal composition is naturally defined
as $f\ten g:P\ten_{H}Q\ra P'\ten_{H}Q'$. We will denote by
$\EtGPD$ the full subbicategory of $\GPD$ consisting of \'{e}tale
Lie groupoids.
\end{ex}

Now let $\C$ and $\D$ be two bicategories. A {\em morphism}
$F:\C\ra\D$ consists of
\begin{enumerate}
\item A function $F:\C_{0}\ra\D_{0}$.

\item Functors $F_{AB}:\C(A,B)\ra\D(F(A),F(B))$ for each pair of
objects of $\C$.
\end{enumerate}
A morphism of bicategories is functorial only up to a family of
natural isomorphisms
\[
\phi_{f,g}:F(f)\circ F(g)\ra F(f\circ
g)\hspace{10pt}\text{and}\hspace{10pt}\phi_{A}:1_{F(A)}\ra
F(1_{A})
\]
for each pair of composable $1$-cells of $\C$ and for each object
of $\C$. These natural isomorphisms have to satisfy some further
natural coherence axioms. Morphisms of bicategories are also
referred to as Lax functors in the literature. Contavariant
versions of morphisms between bicategories can be defined
analogously.

Our main examples of morphisms between bicategories come from the
correspondence between representations of a Lie groupoid and
modules of finite type and of constant rank over its convolution
algebra.

We will first describe the contravariant morphism
\[
\Rep:\EtGPD\ra\AdCat
\]
from the Morita bicategory of \'{e}tale Lie groupoids to the
2-category of additive categories. To any \'{e}tale Lie groupoid
$G$ we assign the additive category $\Rep(G)$ of representations
of $G$. If $H$ is another \'{e}tale Lie groupoid and if $P$ is a
principal $H$-bundle over $G$, then the functor
$\Rep(P):\Rep(H)\ra\Rep(G)$ is defined by
\[
\Rep(P)(E)=P\ten_{H}E\hspace{10pt}\text{and}\hspace{10pt}\Rep(P)(\phi):P\ten_{H}E\ra
P\ten_{H}F
\]
for any $E\in\Rep(H)$ respectively any $\phi\in\Rep(H)(E,F)$. The
morphism $\Rep(P)(\phi)$ of representations of $G$ is explicitely
defined by $\Rep(P)(\phi)(p\ten v)=p\ten\phi(v)$. Furthermore, if
$f:P\ra P'$ is an isomorphism of principal $H$-bundles over $G$,
we define a natural transformation
$\Rep(f):\Rep(P)\Rightarrow\Rep(P')$ by assigning the morphism
\[
\Rep(f)_{E}=f\ten\id:P\ten_{H}E\ra P'\ten_{H}E
\]
of representations of the groupoid $G$ to the representation
$E\in\Rep(H)$. So defined morphism of bicategories is functorial
up to the family of natural isomorphisms
$\phi_{P,Q}:\Rep(P)\circ\Rep(Q)\ra\Rep(P\ten_{H}Q)$ defined by the
natural maps
\[
\phi_{P,Q}(E):P\ten_{H}(Q\ten_{K}E)\ra(P\ten_{H}Q)\ten_{K}E
\]
for any $E\in\Rep(K)$. Natural transformations for the identity
$1$-cells can be defined in a similar fashion by applying the
representation maps.

Our second example is the contravariant morphism
\[
\Mod:\EtGPD\ra\AdCat
\]
from the Morita bicategory of \'{e}tale Lie groupoids to the
$2$-category of additive categories. For each \'{e}tale Lie
groupoid $G$ let $\Mod(G)$ denote the category of modules over the
convolution algebra $\Cc(G)$ that are of finite type and of
constant rank. If $P$ is a principal $H$-bundle over $G$ the
functor $\Mod(P):\Mod(H)\ra\Mod(G)$ is simply the restriction of
the functor of tensoring by $\Cc(P)$
\[
\Cc(P)\ten_{\Cc(H)}-:{}_{\phantom{}H}\mathsf{Mod}\ra{}_{\phantom{}G}\mathsf{Mod}
\]
Tensoring by the right $\Cc(H)$-module $\Cc(P)$ maps the category
of left $\Cc(H)$-modules to the category of left $\Cc(G)$-modules
since $\Cc(P)$ is a $\Cc(G)$-$\Cc(H)$-bimodule. If $\cM\in\Mod(H)$
is a $\Cc(H)$-module of finite type and of constant rank, then
$\Cc(P)\ten_{\Cc(H)}\cM$ is a $\Cc(G)$-module of finite type and
of constant rank by Proposition \ref{Lemma_natural transformation}
which shows that $\Mod(P)$ is a well defined functor. Natural
isomorphisms $\phi_{P,Q}:\Mod(P)\circ\Mod(Q)\ra\Mod(P\ten_{H}Q)$
are defined by the maps
\[
\Omega_{P,Q}\ten\id:\Cc(P)\ten(\Cc(Q)\ten\Gc(E))\ra\Cc(P\ten_{H}Q)\ten\Gc(E),
\]
where $\Omega_{P,Q}:\Cc(P)\ten_{\Cc(H)}\Cc(Q)\ra \Cc(P\ten_{H}Q)$
is the isomorphism of $\Cc(G)$-$\Cc(K)$-bimodules as defined in
\cite{Mrc99}, see also \cite{KalMrc}. The identity coherence
isomorphisms are provided by the actions of convolution algebras
on the modules.

It is straightforward to check that so defined families of
functors and natural isomorphisms satisfy the axioms for the
morphisms between bicategories.

\subsection{Natural equivalence of the morphisms Rep and Mod}

We will show in this section how one can interpret the
Serre-Swan's correspondence as a natural transformation between
the contravariant morphisms $\Rep$ and $\Mod$ from the Morita
bicategory of \'{e}tale Lie groupoids to the $2$-category of
additive categories.

To this effect we first review the definition of a natural
transformation between two morphisms of bicategories
\cite{Lein,Lein2}. Let $\C$ and $\D$ be two bicategories and let
$F,G:\C\ra\D$ be morphisms between them. A {\em natural
transformation} $\sigma:F\ra G$ consists of
\begin{enumerate}
\item For each $A\in\C_{0}$ a $1$-cell $\sigma_{A}:F(A)\ra G(A)$.

\item For each $1$-cell $f:A\ra B$ in $\C$ a $2$-cell
$\sigma_{f}:G(f)\circ\sigma_{A}\Rightarrow\sigma_{B}\circ F(f)$
\[
\begin{CD}
F(A) @>F(f)>>F(B)  \\
@V\sigma_{A}VV @VV\sigma_{B}V \\
G(A) @>G(f)>>G(B)
\end{CD}
\]
\end{enumerate}
such that $\sigma_{f}$ is natural in $f$ and satisfies coherence
axioms as in \cite{Lein2}. We begin with a proposition that is of
some interest independently of our discussion.

\begin{prop}\label{Lemma_natural transformation}
Let $G$ and $H$ be \'{e}tale Lie groupoids and let $P$ be a
principal $H$-bundle over $G$. For any representation $E$ of $H$
there exists a natural isomorphism
\[
\sigma_{P}(E):\Cc(P)\ten_{\Cc(H)}\Gc(E)\ra\Gc(P\ten_{H}E)
\]
of $\Cc(G)$-modules.
\end{prop}
\begin{proof}
We first define a bilinear map
\[
\sigma_{P}(E):\Cc(P)\times\Gc(E)\ra\Gc(P\ten_{H}E)
\]
by the formula
\[
(\sigma_{P}(E)(f,u))(x)=\sum_{\pi(p)=x}f(p)(p\ten u(\phi(p)))
\]
for $f\in\Cc(P)$ and $u\in\Gc(E)$. It is not too hard to check
that the map $\sigma_{P}(E)$ is well defined and that it induces a
homomorphism
\[
\sigma_{P}(E):\Cc(P)\ten_{\Cc(H)}\Gc(E)\ra\Gc(P\ten_{H}E)
\]
of $\Cc(G)$-modules, which we claim to be an isomorphism. It
suffices to show that $\sigma_{P}(E)$ is a bijective map.

We first consider the case when the bundle $P$ is trivial, i.e.
$P=P(\psi)$ for some smooth functor $\psi:G\ra H$. The
representation $P\ten_{H}E$ of $G$ is then isomorphic to the
representation $\psi^{\ast}E$ of $G$ via the isomorphism
$f:P\ten_{H}E\ra\psi^{\ast}E$. In particular, $f$ is an
isomorphism of vector bundles over $M$, so we obtain the
isomorphism
\[
\Gc(f):\Gc(P\ten_{H}E)\ra\Gc(\psi_{0}^{\ast}E)
\]
of $\Cc(M)$-modules. Now consider the representation $E$ of $H$ as
a vector bundle over $N$. It is a classical result (see
\cite{GrHaVa}) that the map
\[
\sigma_{\psi_{0}}(E):\Cc(M)\ten_{\Cc(N)}\Gc(E)\ra\Gc(\psi_{0}^{\ast}E),
\]
defined analogously as the map $\sigma_{P}(E)$, is an isomorphism
of $\Cc(M)$-modules. Finally, since $P=M\times_{N}^{\psi,t}H$ is a
trivial bundle, we have by \cite{Mrc99} the isomorphism
$\Omega_{M,H}:\Cc(P)\cong\Cc(M)\ten_{\Cc(N)}\Cc(H)$ of
$\Cc(M)$-$\Cc(H)$-bimodules which gives us an isomorphism
\[
\iota:\Cc(P)\ten_{\Cc(H)}\Gc(E)\ra\Cc(M)\ten_{\Cc(N)}\Gc(E)
\]
of $\Cc(M)$-modules. We can collect all these isomorphisms into
the following commutative diagram of homomorphisms of
$\Cc(M)$-modules
\[
\begin{CD}
\Cc(P)\ten_{\Cc(H)}\Gc(E) @>\sigma_{P}(E)>>\Gc(P\ten_{H}E)  \\
@V\iota V\cong V @V\cong V\Gc(f)V \\
\Cc(M)\ten_{\Cc(N)}\Gc(E) @>\sigma_{\psi_{0}}(E)>\cong
>\Gc(\psi_{0}^{\ast}E)
\end{CD}
\]
Since the remaining three maps are bijective, the map
$\sigma_{P}(E)$ is bijective as well.

A principal $H$-bundle $P$ over $G$ is in general only locally
trivial \cite{MoMr2}. Let $U$ be an open subset of $M$ such that
$P|_{U}$ is a trivial $H$-bundle. We then have a natural injective
homomorphism $i:\Cc(P|_{U})\ra\Cc(P)$ of right $\Cc(H)$-modules
which induces an injective homomorphism
\[
i\ten\id:\Cc(P|_{U})\ten_{\Cc(H)}\Gc(E)\ra\Cc(P)\ten_{\Cc(H)}\Gc(E)
\]
of abelian groups. Injectivity of $i\ten\id$ follows from the fact
that $\Cc(P)$ is a locally unital $\Cc(M)$-module. Namely, for any
$w=\sum f_{i}\ten u_{i}\in\Cc(P|_{U})\ten_{\Cc(H)}\Gc(E)$ there
exists a function $f\in\Cc(M)$ with support in $U$ such that
$ff_{i}=f_{i}$ for each $i$. Left action
$\mu_{f}:\Cc(P)\ten_{\Cc(H)}\Gc(E)\ra\Cc(P|_{U})\ten_{\Cc(H)}\Gc(E)$
by $f$ is then a homomorphism of abelian groups such that
$\mu_{f}((i\ten\id)(w))=w$. We now have the following commuting
square of homomorphisms of abelian groups
\[
\begin{CD}
\Cc(P|_{U})\ten_{\Cc(H)}\Gc(E) @>\sigma_{P|_{U}}(E)>\cong>\Gc(P|_{U}\ten_{H}E)  \\
@Vi\ten\id VV @VVV \\
\Cc(P)\ten_{\Cc(H)}\Gc(E) @>\sigma_{P}(E)>>\Gc(P\ten_{H}E)
\end{CD}
\]
where the right vertical map is defined by trivially extending the
sections outside of the set $U$. Note that both the vertical maps
are injective. The injectivity of the map $\sigma_{P}(E)$ will now
follow from the injectivity of the map $\sigma_{P|_{U}}(E)$.
Indeed, suppose that $\eta\in\Cc(P)\ten_{\Cc(H)}\Gc(E)$ is such
that $\sigma_{P}(E)(\eta)=0$. Since the module $\Cc(P)$ is locally
unital, there exists a function $f\in\Cc(M)$ such that
$f\eta=\eta$ and a decomposition $f=\sum f_{i}$ into functions
with supports contained in open sets $U_{i}$ such that the bundle
$P|_{U_{i}}$ is trivial for each $i$. It now follows from the
$\Cc(M)$-linearity of the map $\sigma_{P}(E)$ and from the
preceding diagram that
\[
\sigma_{P|_{U_{i}}}(E)(f_{i}\eta)=\sigma_{P}(E)(f_{i}\eta)=f_{i}\sigma_{P}(E)(\eta)=0.
\]
The injectivity of the maps $\sigma_{P|_{U_{i}}}(E)$ now implies
that $f_{i}\eta=0$ for each $i$ and thus $\eta=0$, which proves
that $\sigma_{P}(E)$ is injective. Surjectivity can be proven by
using similar arguments and by noting that the sets
$\Cc(P|_{U})\ten_{\Cc(H)}\Gc(E)$, as $U$ varies over some
trivializing cover of the bundle $P$, generate the $\Cc(G)$-module
$\Cc(P)\ten_{\Cc(H)}\Gc(E)$.
\end{proof}

Let us now return to our morphisms $\Rep$ and $\Mod$ of
bicategories. Define a natural transformation
$\sigma:\Rep\Rightarrow\Mod$ by defining a functor
\[
\sigma_{G}=(\Gc)_{G}
\]
for each \'{e}tale Lie groupoid $G$ and a natural transformation
\[
\sigma_{P}:\Mod(P)\circ\sigma_{H}\Rightarrow\sigma_{G}\circ\Rep(P)
\]
for each $1$-cell $P\in\EtGPD(G,H)$ of the bicategory $\EtGPD$.
Taking into account Proposition \ref{Lemma_natural transformation}
it is now straightforward to verify the following theorem.

\begin{theo}\label{Natural equivalence of morphisms}
The morphisms $\Rep$ and $\Mod$ from the Morita bicategory of Lie
groupoids to the bicategory of additive categories are naturally
equivalent. Natural equivalence is given by the the family of
functors $(\Gc)_{G}:\Rep(G)\ra\Mod(G)$ for any groupoid $G$ and
the family of transformations
$\sigma_{P}:\Mod(P)\circ\sigma_{H}\Rightarrow\sigma_{G}\circ\Rep(P)$
for any principal $H$-bundle $P$ over $G$.
\end{theo}

\noindent {\bf Acknowledgements.} I would like to thank J.
Mr\v{c}un for many helpful discussions and advice related to the
paper.

\end{document}